\documentclass[12pt]{article}

      \usepackage{latexsym}
         \usepackage[reqno, namelimits, sumlimits]{amsmath}
         \usepackage{amssymb, amsfonts}
         \usepackage{amsthm}

 \newtheorem{theorem}{Theorem}[section]

 \newtheorem{lemma}[theorem]{Lemma}

 \newtheorem{pro}[theorem]{Proposition}
 
\title{ A necessary condition of possible blowup for the Navier-Stokes system in half-space.}

 \author{T. Barker, G. Seregin }

\begin{document}
\maketitle

\setcounter{equation}{0}
\section{Introduction}

The question that is addressed in the paper is as follows. Let us consider the initial boundary value problem for the Navier-Stokes system in the space-time domain $Q_+=\Omega\times ]0,\infty[$ for vector-valued function $v=(v_1,v_2,v_3)=(v_i)$ and scalar function $q$, satisfying the equations
\begin{equation}\label{directsystem}
\partial_tv+v\cdot\nabla v-\Delta v=-\nabla q,\qquad\mbox{div}\,v=0
\end{equation}
in $Q_+$,
the boundary conditions
\begin{equation}\label{directbc}
v=0\end{equation}
on $\Omega\times [o,\infty[$,
and the initial conditions
\begin{equation}\label{directic}
v(\cdot,0)=v_0(\cdot)\end{equation}
in $\Omega$. It is assumed that the initial velocity field $v_0$ is smooth, compactly supported, and divergence free in $\Omega$, i.e., $v_0$ belongs to the space $C^\infty_{0,0}(\Omega)$, and that $\Omega$ is a domain in $\mathbb R^3$ with sufficiently smooth boundary. Our main aim is to study whether or not the velocity field $v$ blows up in a finite time, in other words, whether or not there exists a finite time $T>0$  such that
\begin{equation}\label{definitionblowup}
\lim\limits_{t\uparrow T}\|v(\cdot,t)\|_{\infty,\Omega}=\infty.\end{equation}
There is a huge number of papers dedicated to this problem. Among them the most relevant to us are the following papers. In the first place, one should mention the classical Leray necessary conditions for $T$ to be a blowup time:
\begin{equation}\label{Lerayrates}
\|v(\cdot,t)\|_{s,\Omega}\geq \frac {c_s}{(T-t)^{\frac {s-3}{2s}}}
\end{equation}
for any $0<t<T$, for all $s>3$, and  for a positive constant $c_s$ depending only on $s$.
Estimates (\ref{Lerayrates}) have been proven by J. Leray in \cite{Le} for $\Omega=\mathbb R^3$ and then by Y. Giga in \cite{Giga1986} for a wide class of domains $\Omega$ including a half space and bounded domains with sufficiently smooth boundaries. However, there is an interesting  marginal case $s=3$, in which no estimate of type (\ref{Lerayrates}) is  known. In papers \cite{ESS2003},\cite{MSh2006} and \cite{S2005}, it has been shown that
\begin{equation}\label{limsup}
\limsup\limits_{t\uparrow T}\|u(\cdot,t)\|_{3,\Omega}=\infty
\end{equation}
for $\Omega=\mathbb R^3$, $\Omega=\mathbb R^3_+:=\{x=(x_i)\in \mathbb R^3:\,\,x_3>0\}$ and for $\Omega$ being a bounded domain with sufficently smooth boundary. Later on, in series of papers \cite{S6}
--\cite{Ser2012}, necessary condition (\ref{limsup}) has been improved for $\Omega=\mathbb R^3$ in the following sense
\begin{equation}\label{lim}
\lim\limits_{t\uparrow T}\|u(\cdot,t)\|_{3,\Omega}=\infty
\end{equation}
The aim of the paper is to prove
\begin{theorem}\label{limhalf}
Necessary condition (\ref{lim}) remains to be true for $\Omega=\mathbb R^3_+$.
\end{theorem}
We also believe that necessary condition (\ref{lim}) holds when $\Omega$ is a bounded domain with sufficiently smooth boundary. The proof of this will be published elsewhere.

We would like to empathise that to prove Theorem \ref{limhalf}, a different approach is used to that of the whole space. Though we focus on the half space, this also provides an alternative to the proof given for the whole space in \cite{Ser2012}.

 The main difficulty in attempts to prove Theorem \ref{limhalf} is as follows. The proof of this statement in the case $\Omega=\mathbb R^3$ consists of two big parts: rescaling,  leading to a certain class of ancient solutions to the Navier-Stokes equations, and a Liouville type theorem for those solutions based on the backward uniqueness. The second part at least conceptually works  in the case of a half space $\mathbb R^3_+$ as well while the first one does not. The reason is that the rescaling and the limiting procedure in the case of the whole space $\mathbb R^3$ give the special type of the so-called local energy ancient solutions to the Navier-Stokes that coincide with Lemarie-Rieusset solutions to the Cauchy problem for the Navier-Stokes equations on some finite time interval. Those solutions have been introduced by Lemarie-Rieusset in \cite{LR1}, see also for some definitions in \cite{KS}. Unfortunately, an analog of Lemarie-Rieusset solutions for a half space is not known yet. In fact, this is an interesting open problem. In this paper, we are able to work without Lemarie-Rieusset type solutions in half space to get a local energy ancient solution to which a Liouville type theorem based on backward uniqueness is applicable.

\setcounter{equation}{0}
\section{A priori estimates}
 Let us consider a sufficiently smooth solution $u$ and $p$ to the Navier-Stokes system  in the space-time strip $Q^+_{-2,0}=\mathbb R^3_+\times ]-2,0[$ to the following initial boundary value problem:
 \begin{equation}\label{system}
\partial_tu+\mbox{div}\,u\otimes u-\Delta u=-\nabla p,\qquad \mbox{div}\,u=0 \end{equation}
in $Q^+_{-2,0}$,
\begin{equation}\label{hombc}u(x',0,t)=0\end{equation}
for $(x',t)\in\mathbb R^2\times [-2,0]$,
$$u(\cdot,-2)=u_0(\cdot)\in L_3(\mathbb R^3_+).$$
We may split the solution into two parts
$$u=u^1+u^2,$$
where $u^1$ and $p^1$ solve the linear problem
$$\partial_tu^1-\Delta u^1=-\nabla p^1,\qquad \mbox{div}\,u^1=0$$
in $Q^+_{-2,0}$,
$$u^1(x',0,t)=0$$
for $(x',t)\in\mathbb R^2\times [-2,0]$,
$$u^1(\cdot,-2)=u_0(\cdot)\in L_3(\mathbb R^3_+).$$
Suppose that
\begin{equation}\label{estforinitialdata}
\|u_0\|_{3,\mathbb R^3_+}\leq M.
\end{equation}

Using Solonnikov estimates for the Green function in a half space, see \cite{{Sol1973}} and \cite{Sol2003UMN}, one can check  all assumptions in the Lemma  of \cite{Giga1986} and state that the following two estimates for $u^1$ are valid:
$$\|u^1\|_{3,\infty,Q^+_{-2,0}}+\|u^1\|_{5,Q^+_{-2,0}}\leq c\|u_0\|_{3,\mathbb R^3_+}\leq cM.$$
Hence, simply by the interpolation, we have
\begin{equation}\label{globalu1}
\|u^1\|_{s,Q^+_{-2,0}} \leq c(s)\|u_0\|_{3,\mathbb R^3_+}\leq c(s)M\end{equation}
for any $s\in [3,5]$. In addition, the above mentioned Solonnikov estimates implies the following inequality
\begin{equation}\label{globgradu1}
\|\nabla u^1(\cdot,t)\|_{3,\mathbb R^3_+}\leq \frac c{\sqrt{t+2}}\|u_0(\cdot)\|_{3,\mathbb R^3_+}\leq  \frac {cM}{\sqrt{t+2}}\end{equation}
for any $t\in ]-2,0[$.

The second counterpart of $u$ satisfies the  non-linear system
$$\partial_tu^2+\mbox{div}\,u\otimes u-\Delta u^2=-\nabla p^2,\qquad \mbox{div}\,u^2=0$$
in $Q^+_{-2,0}$, the boundary conditions 
$$u^2(x',0,t)=0$$
for $(x',t)\in\mathbb R^2\times [-2,0]$, and the initial conditions
$$u^2(\cdot,-2)=0$$
in $R^3_+$.

The standard energy approach to the second system gives
$$ \partial_t\|u^2(\cdot,t)\|^2_{2,\mathbb R^3_+}+2\|\nabla u^2(\cdot,t)\|^2_{2,\mathbb R^3_+}=$$
$$= 2\int\limits_{\mathbb R^3_+}u\otimes u:\nabla u^2dx ds=I_1+I_2+I_3+I_4,$$
where
$$I_1=2\int\limits_{\mathbb R^3_+}u^1\otimes u^1:\nabla u^2dx,$$
$$I_2=2\int\limits_{\mathbb R^3_+}u^1\otimes u^2:\nabla u^2dx,$$
$$I_3=2\int\limits_{\mathbb R^3_+}u^2\otimes u^1:\nabla u^2dx=0,$$
$$I_4=2\int\limits_{\mathbb R^3_+}u^2\otimes u^2:\nabla u^2dx=0.$$
Next, let us  consequently evaluate terms on the right hand side of the energy identity. For the first term, we have
$$|I_1|\leq c\|u^1(\cdot,t)\|^2_{4,R^3_+}\|\nabla u^2(\cdot,t)\|_{2,\mathbb R^3_+}.$$
The second term can be treated as follows:
$$|I_2|\leq c\|u^1(\cdot,t)\otimes u^2(\cdot,t)\|_{2,\mathbb R^3_+}\|\nabla u^2(\cdot,t)\|_{2,\mathbb R^3_+}\leq$$
$$\leq c\|u^1(\cdot,t)\|_{5,\mathbb R^3_+}\|u^2(\cdot,t)\|_{\frac {10}3,\mathbb R^3_+}\|\nabla u^2(\cdot,t)\|_{2,\mathbb R^3_+}.
$$
Applying the known multiplicative inequality to the second factor in the right hand side of the latter bound, we find
$$|I_2|\leq c\|u^1(\cdot,t)\|_{5,\mathbb R^3_+}\|u^2(\cdot,t)\|^{\frac 25}_{2,\mathbb R^3_+}\|\nabla u^2(\cdot,t)\|^\frac 85_{2,\mathbb R^3_+}.$$
Letting
$$y(t):=\|u^2(\cdot,t)\|^2_{2,\mathbb R^3_+}$$
and using the Young inequality, we find
$$y'(t)+\|\nabla u^2(\cdot,t)\|^2_{2,\mathbb R^3_+}\leq
c\|u^1(\cdot,t)\|^5_{5,\mathbb R^3_+}y(t)+c\|u^1(\cdot,t)\|^4_{4,R^3_+}.$$
Next, elementary arguments lead to the inequality
$$\Big(y(t)\exp{\Big(-\int\limits^t_{-2}}\|u^1(\cdot,s)\|^5_{5,\mathbb R^3_+}ds\Big)\Big)'\leq c\exp{\Big(-\int\limits^t_{-2}}\|u^1(\cdot,s)\|^5_{5,\mathbb R^3_+}ds\Big)\|u^1(\cdot,t)\|^4_{4,R^3_+}.$$
So,
$$y(t)\leq c\int\limits^t_{-2}\exp{\Big(\int\limits^t_\tau}\|u^1(\cdot,s)\|^5_{5,\mathbb R^3_+}ds\Big)\|u^1(\cdot,\tau)\|^4_{4,R^3_+}d\tau\leq$$
$$\leq c\|u_0\|^4_{3,\mathbb R^3_+}\exp\Big({c}\|u_0\|^5_{3,\mathbb R^3_+}\Big)$$
and
\begin{equation}\label{energyu2}\|\nabla u^2\|^2_{2,Q^+_{-2,0}}\leq c\|u_0\|^4_{3,\mathbb R^3_+}+c\|u_0\|^9_{3,\mathbb R^3_+}\exp\Big({c}\|u_0\|^5_{3,\mathbb R^3_+}\Big)\end{equation}
From these estimates and from the multiplicative inequality, see also (\ref{estforinitialdata}),
one can deduce that
\begin{equation}\label{globalu2}
\|u^2\|_{s,Q^+_{-2,0}}\leq C(s,M)\end{equation}with any $s\in [2,\frac {10}3]$.

Let us fix a smooth cut-off function $\chi(t)$ so that $\chi(t)=1$ if $-3/2<t<1$ and $\chi(t)=0$ if $-2<t<-7/4$.
Then, we may split $\chi u^2$ and $\chi p^2$ in the following way:
$$\chi u^2=u^{2,1}+u^{2,2}+u^{2,3}+u^{2,4} +u^{2,5}$$
and
$$\chi p^2=p^{2,1}+p^{2,2}+p^{2,3}+p^{2,4}+p^{2,5}$$
so that, for $i=1,2,3,4,5$,
$$\partial_tu^{2,i}-\Delta u^{2,i}+\nabla p^{2,i}= f^i,\qquad \mbox{div}\,u^{2,i}=0$$
in $Q^+_{-2,0}$,
$$u^{2,i}(x',0,t)=0$$
for all $(x',t)\in \mathbb R^2\times [-2,0]$ and
$$u^{2,i}(x,0)=0$$
for $x\in \mathbb R^2_+$, where
$$f^1:=\chi'u^2,\quad f^2:=-\chi u^2\cdot\nabla u^2,\quad f^3:=-\chi u^2\cdot \nabla u^1, $$$$ f^4:=-\chi u^1\cdot\nabla u^2,\quad f^5:=-\chi u^1\cdot \nabla u^1.$$

We start with evaluation of $u^{2,1}$. Our main tool here is the Solonnikov coercive estimates of the linear theory. In particular, it  follows from (\ref{globalu2}) that
\begin{equation}\label{timeu21}
\|\partial_tu^{2,1}\|_{s,Q^+_{-2,0}}+\|\nabla^2u^{2,1}\|_{s,Q^+_{-2,0}}+\|\nabla p^{2,1}\|_{s,Q^+_{-2,0}}\leq $$$$\leq c(s)\|f^1\|_{s,Q^+_{-2,0}}\leq c(s)\|u^2\|_{s,Q^+_{-2,0}} \leq C(s,
M)\end{equation}
for any $s\in [2,10/3]$. To estimate the second counter-part $u^{2,2}$, one can use the standard consequence of the energy bounds and find
\begin{equation}\label{time22}
\|\partial_tu^{2,2}\|_{s,l,Q^+_{-2,0}}+\|\nabla^2u^{2,2}\|_{s,l,Q^+_{-2,0}}+\|\nabla p^{2,2}\|_{s,l,Q^+_{-2,0}}\leq $$$$\leq c(s,l)\|f^2\|_{s,l,Q^+_{-2,0}}\leq C(s,
M)\end{equation}
provided that
$$\frac 3s+\frac 2l=4.$$

Next, for $i=3,5$, it follows from (\ref{globalu1}), (\ref{globgradu1}), and (\ref{globalu2}) that
\begin{equation}\label{time23and25}
\max\limits_{i=3,5}\Big(\|\partial_tu^{2,i}\|_{3/2,Q^+_{-2,0}}+\|\nabla^2u^{2,i}\|_{3/2,Q^+_{-2,0}}+\|\nabla p^{2,i}\|_{3/2,Q^+_{-2,0}}\Big)\leq $$$$\leq c(\|f^3\|_{3/2,Q^+_{-2,0}}+\|f^5\|_{3/2,Q^+_{-2,0}})\leq C(
M).\end{equation}
Finally, applying H\"older inequality, we have a bound for $u^{2,4}$:
\begin{equation}\label{timeu24}
\|\partial_tu^{2,4}\|_{6/5,3/2,Q^+_{-2,0}}+\|\nabla^2u^{2,4}\|_{6/5,3/2,Q^+_{-2,0}}+\|\nabla p^{2,4}\|_{6/5,3/2,Q^+_{-2,0}}\leq $$$$\leq c\|f^4\|_{6/5,3/2,Q^+_{-2,0}}\leq c2^\frac 16\|u^1\|_{3,\infty,Q^+_{-2,0}}\|\nabla u^2\|_{2,Q^+_{-2,0}}\leq C(
M)\end{equation}

As to $u^1$, we let $v^1:=\chi u^1$ and $q^1:=\chi p^1$ and find
$$\partial_t v^1-\Delta v^1+\nabla q^1=\chi'u^1,\qquad \mbox{div}\,v^1=0$$
in $Q^+_{-2,0}$,
$$v^1(x',0,t)=0$$
for all $(x',t)\in \mathbb R^2\times ]-2,0[$ and
$$v^1(\cdot,0)=0$$
for all $x\in \mathbb R^3_+$. The same arguments as above lead to the estimate
\begin{equation}\label{timeu1}\|\partial_tv^{1}\|_{3,Q^+_{-2,0}}+\|\nabla^2v^{1}\|_{3,Q^+_{-2,0}}+\|\nabla q^{1}\|_{3,Q^+_{-2,0}}\leq $$$$\leq c\|u^1\|_{3,Q^+_{-2,0}}\leq C(
M).\end{equation}


In what follows, we are going to use the following Poincare type inequalities:
\begin{equation}\label{pres3and5}
\int\limits^0_{-3/2}\int\limits_{B(x_0,R)}|p^{2,i}-[p^{2,i}]_{B(x_0,R)}|^\frac 32 dx dt\leq cR^\frac 32
\int\limits^0_{-3/2}\int\limits_{B(x_0,R)}|\nabla p^{2,i}|^\frac 32dx dt\end{equation}
for $i=3,5$;
\begin{equation}\label{pres1}\int\limits^0_{-3/2}\int\limits_{B(x_0,R)}|p^{2,1}-[p^{2,1}]_{B(x_0,R)}|^\frac 32dx dt\leq c R^\frac {9}4
\Big(\int\limits^0_{-3/2}\int\limits_{B(x_0,R)}|\nabla p^{2,1}|^2dx dt\Big)^\frac 34;\end{equation}
\begin{equation}\label{pres2}
\int\limits^0_{-3/2}\int\limits_{B(x_0,R)}|p^{2,2}-[p^{2,2}]_{B(x_0,R)}|^\frac 32dx dt\leq c R^\frac 12\int\limits^{0}_{-3/2}\Big(\int\limits_{B(x_0,R)}|\nabla p^{2,2}|^\frac 98dx\Big)^\frac 43dt;\end{equation}
\begin{equation}\label{pres4}\int\limits^0_{-3/2}\int\limits_{B(x_0,R)}|p^{2,4}-[p^{2,4}]_{B(x_0,R)}|^\frac 32dx dt\leq c R^\frac 34\int\limits^{0}_{-3/2}\Big(\int\limits_{B(x_0,R)}|\nabla p^{2,4}|^\frac 65dx\Big)^\frac 54dt;\end{equation}
\begin{equation}\label{pres01}\int\limits^0_{-3/2}\int\limits_{B(x_0,R)}|p^{1}-[p^{1}]_{B(x_0,R)}|^\frac 32dx dt\leq c R^2\Big(\int\limits^{0}_{-3/2}\int\limits_{B(x_0,R)}|\nabla p^{1}|^3dxdt\Big)^\frac 12.\end{equation}
All the formulae are valid provided $B(x_0,R)\in \mathbb R^3_+$. They are also valid if we replace $B(x_0,R)$
with semi-balls $B^+(x_0,R)$ assuming that $x_{03}=0$.

\setcounter{equation}{0}
\section{Passage to the limit}

Suppose that
we have a sequence of sufficiently smooth functions $u^{(n)}$ and $p^{(n)}$ defined in the domain $Q^+_{-2,0}=\mathbb R^3_+\times ]-2,0[$ that are solutions to the following initial boundary value problem:
\begin{equation}\label{nsystem}
\partial_tu^{(n)}+\mbox{div}\,u^{(n)}\otimes u^{(n)}-\Delta u^{(n)}=-\nabla p^{(n)},\qquad \mbox{div}\,u^{(n)}=0
\end{equation}
in $Q^+_{-2,0}$,
\begin{equation}\label{nbc} u^{(n)}(x',0,t)=0\end{equation}
for $(x',t)\in\mathbb R^2\times [-2,0]$,
\begin{equation}\label{nic}
u^{(n)}(\cdot,-2)=u^{(n)}_0(\cdot)\in L_3(\mathbb R^3_+).\end{equation}
It is supposed also that
$$u^{(n)}_0\rightharpoonup u_0$$
in $L_3(\mathbb R^3_+)$.

We let
$$M:=\sup\limits_n\|u^{(n)}_0\|_{3,\mathbb R^3_+}<\infty.$$
\begin{pro}\label{convergence}
There exist subsequences still denoted in the same way with the following properties:
\begin{equation}\label{weakconvergence}
u^{(n)}\rightharpoonup u
\end{equation}
in $L_\frac {10}3(Q^+_{-2,0})$,
\begin{equation}\label{weakconvergencenabla}\nabla u^{(n)}\rightharpoonup \nabla u\end{equation}
in $L_2(B^+(R)\times ]-2+\delta,0[)$ for any $R>0$ and any $0<\delta<2$,
\begin{equation}\label{strongconvergence}
u^{(n)}\rightarrow u
\end{equation}
in $L_3(B^+(R)\times ]-3/2,0[)$ for any $R>0$  and
\begin{equation}\label{weakconverpressure}
p^{(n)}\rightharpoonup p
\end{equation}
in $L_{\frac 32}(B^+(R)\times ]-3/2,0[)$ for any $R>0$. 

Functions $u$ and $p$ satisfy (\ref{system}) in $\mathbb R^3_+\times ]-3/2,0[$ and (\ref{hombc}) for $(x',t)\in \mathbb R^2\times ]-3/2,0[$.

For the pressure $p$, the following global estimates are valid:
$$p=p^1+p^2$$
and
$$p^{2}=\sum\limits^5_{i=1}p^{2,i}.
$$
with the estimates
$$\|\nabla p^1\|_{3,\mathbb R^3_+\times ]-3/2,0[}+\|\nabla p^{2,1}\|_{2,\mathbb R^3_+\times ]-3/2,0[}+$$
\begin{equation}\label{gradpressure}
+\|\nabla p^{2,2}\|_{9/8,3/2,\mathbb R^3_+\times ]-3/2,0[}
+\|\nabla p^{2,3}\|_{3/2,\mathbb R^3_+\times ]-3/2,0[}+\end{equation}
$$
+\|\nabla p^{2,4}\|_{6/5,3/2,\mathbb R^3_+\times ]-3/2,0[}
+\|\nabla p^{2,5}\|_{3/2,\mathbb R^3_+\times ]-3/2,0[}<\infty.$$

Moreover, for any $R>0$, the limits pair $u$ and $p$ satisfies the local energy inequality
$$\int\limits_{B(x_0,R)\cap\mathbb R^3_+}|\varphi^2(x,t)|u(x,t)|^2dx +2\int\limits^{t}_{t_0-R^2}\int\limits_{B(x_0,R)\cap\mathbb R^3_+}\varphi^2|\nabla u|^2dxdt   \leq $$
\begin{equation}\label{locenergy}
\leq\int\limits^{t}_{t_0-R^2}\int\limits_{B(x_0,R))\cap\mathbb R^3_+}\Big(|u|^2(\partial_t\varphi^2+\Delta \varphi^2u^\cdot \varphi^2(|u|^2+2p)\Big) dxds
\end{equation}
for all $-3/2<t_0-R^2<t\leq t_0\leq0$, for all $x_0\in R^3$, and for all $\varphi\in C^\infty_0(B(x_0,R)\times ]t_0-R^2,t_0+R^2[)$.\end{pro}
\textbf{Proof} Obviously,  we may assume, without loss of generality,  that (\ref{weakconvergence}) follows from (\ref{globalu1}) and (\ref{globalu2}). Moreover,
the limit function obeys the estimate
\begin{equation}\label{globallimitu}
\|u\|_{10/3,Q^+_{-2,0}}<\infty.\end{equation}
Obviously,  (\ref{globgradu1}) and (\ref{energyu2}) imply (\ref{weakconvergencenabla}). 
From 
(\ref{timeu21})--
(\ref{timeu1}), we can deduce (\ref{strongconvergence}).

Now, let us treat the pressure $p^{(n)}$, using the decomposition of the previous section
$$p^{(n)}=p^{(n)1}+p^{(n)2},$$
where
$$p^{(n)2}=\sum\limits^5_{i=1}p^{(n)2,i}.
$$
Then, using (\ref{timeu21})--(\ref{timeu1}) and (\ref{pres3and5})--(\ref{pres01}), we can justify (\ref{weakconverpressure}) and (\ref{gradpressure}).

Since functions $u^{(n)}$ and $p^{(n)}$ satisfy the local energy inequality, i.e.,
$$\int\limits_{B(x_0,R)\cap\mathbb R^3_+}|\varphi^2(x,t)|u^{(n)}(x,t)|^2dx +2\int\limits^{t}_{t_0-R^2}\int\limits_{B(x_0,R)\cap\mathbb R^3_+}\varphi^2|\nabla u^{(n)}|^2dxdt   \leq $$$$\leq\int\limits^{t}_{t_0-R^2}\int\limits_{B(x_0,R)\cap\mathbb R^3_+}\Big(|u^{(n)}|^2(\partial_t\varphi^2+\Delta \varphi^2)+u^{(n)}\cdot \varphi^2(|u^{(n)}|^2+2p^{(n)})\Big) dxds$$
for all $-2< t_0-R^2<t\leq t_0\leq 0$, for all $x_0\in \mathbb R^3$, and for all $\varphi\in C^\infty_0(B(x_0,R)\times ]t_0-R^2,t_0+R^2[)$, we can find (\ref{locenergy}) by passing to the limits and taking into account (\ref{strongconvergence}) and (\ref{weakconverpressure}).

\begin{pro}\label{decaypro} Let $u$ and $p$ be a limit function from Proposition \ref{convergence}. There exist a number $R_1>0$ such that
\begin{equation}\label{maxu}
|u(x,t)|\leq c
\end{equation}
for all $(x,t)\in (\mathbb R^3_+\setminus B^+(R_1))\times ]-5/4,0[$ and for some universal constant $c$.
Moreover, given $\delta>0$,
\begin{equation}\label{maxgradu}
|\nabla u(x,t)|\leq c_1(\delta)
\end{equation}
for all $(x,t)\in (\mathbb R^3_{+\delta}\setminus B^+(R_1))\times ]-5/4,0[$. Here, $\mathbb R^3_{+\delta}:=\mathbb R^3_+\cap\{x_3>\delta\}$.
\end{pro}
\textbf{Proof} By (\ref{gradpressure}), 
 we can state that
$$\int\limits^0_{-3/2}\int\limits_{\mathbb R^3_+\setminus B^+(R)}|u|^3dxdt+
\Big(\int\limits^{0}_{-3/2}\int\limits_{\mathbb R^3_+\setminus B^+(R)}|\nabla p^{1}|^3dxdt\Big)^\frac 12+
$$$$+\Big(\int\limits^0_{-3/2}\int\limits_{\mathbb R^3_+\setminus B^+(R)}|\nabla p^{2,1}|^2dx dt\Big)^\frac 34+\int\limits^{0}_{-3/2}\Big(\int\limits_{\mathbb R^3_+\setminus B^+(R)}|\nabla p^{2,2}|^\frac 98dx\Big)^\frac 43dt+$$
$$+\sum\limits_{i=1}^2\limits\int\limits^0_{-3/2}\int\limits_{\mathbb R^3_+\setminus B^+(R)}|\nabla p^{2,2i+1}|^\frac 32dx dt+\int\limits^{0}_{-3/2}\Big(\int\limits_{\mathbb R^3_+\setminus B^+(R)}|\nabla p^{2,4}|^\frac 65dx\Big)^\frac 54\to 0$$
as $R\to\infty$.

Given $\varepsilon>0$, there exists a positive number $R_1>0$ such that
$$\frac 1{r^2}\int\limits_{Q(z_0,r)}(|u|^3+|p-[p]_{B(x_0,r)}|^\frac 32)dxdt \leq $$
$$\leq \frac 1{r^2}\int\limits_{Q(z_0,r)}|u|^3dx dt+\frac c{r^2}\int\limits_{Q(z_0,r)}|p^1-[p^1]_{B(x_0,r)}|^\frac 32dxdt+$$
$$+\frac c{r^2}\sum^5_{i=1}\int\limits_{Q(z_0,r)}|p^{2,i}-[p^{2,i}]_{B(x_0,r)}|^\frac 32dxdt\leq \frac 1{r^2}\int\limits_{Q(z_0,r)}|u|^3dx dt+$$
$$+c\Big(\int\limits^{0}_{-3/2}\int\limits_{ B(x_0,r)}|\nabla p^{1}|^3dxdt\Big)^\frac 12+ c r^\frac {1}4
\Big(\int\limits^0_{-3/2}\int\limits_{B(x_0,r)}|\nabla p^{2,1}|^2dx dt\Big)^\frac 34+$$
$$+ cr^{-\frac 32}\int\limits^{0}_{-3/2}\Big(\int\limits_{B(x_0,r)}|\nabla p^{2,2}|^\frac 98dx\Big)^\frac 43dt+ c r^{-\frac 54}\int\limits^{0}_{-3/2}\Big(\int\limits_{B(x_0,r)}|\nabla p^{2,4}|^\frac 65dx\Big)^\frac 54dt+$$$$
+ cr^{-\frac 12}\sum\limits^2_{i=1}
\int\limits^0_{-3/2}\int\limits_{B(x_0,r)}|\nabla p^{2,2i+1}|^\frac 32dx dt<\varepsilon  $$
with $r=1/100$ and $Q(z_0,r)\in (\mathbb R^3_+\setminus B^+(R_1/2))\times ]-3/2,0[$.

The same can be done for boundary points:
$$\frac 1{\varrho^2}\int\limits_{Q^+(z_0,\varrho)}(|u|^3+|p-[p]_{B^+(x_0,\varrho)}|^\frac 32)dxdt
\leq \frac 1{\varrho^2}\int\limits_{Q^+(z_0,\varrho)}|u|^3dx dt+$$
$$+c\Big(\int\limits^{0}_{-3/2}\int\limits_{ B^+(x_0,\varrho)}|\nabla p^{1}|^3dxdt\Big)^\frac 12+ c \varrho^\frac {1}4
\Big(\int\limits^0_{-3/2}\int\limits_{B^+(x_0,\varrho)}|\nabla p^{2,1}|^2dx dt\Big)^\frac 34+$$
$$+ c\varrho^{-\frac 32}\int\limits^{0}_{-3/2}\Big(\int\limits_{B^+(x_0,\varrho)}|\nabla p^{2,2}|^\frac 98dx\Big)^\frac 43dt+ c \varrho^{-\frac 54}\int\limits^{0}_{-3/2}\Big(\int\limits_{B^+(x_0,\varrho)}|\nabla p^{2,4}|^\frac 65dx\Big)^\frac 54dt+$$$$
+ c\varrho^{-\frac 12}\sum\limits^2_{i=1}
\int\limits^0_{-3/2}\int\limits_{B^+(x_0,\varrho)}|\nabla p^{2,2i+1}|^\frac 32dx dt<\varepsilon  $$
with $\varrho=1/10$ and $Q^+(z_0,\varrho):=B^+(x_0,\varrho)\times ]t_0-\varrho^2,t_0[\in (\mathbb R^3_+\setminus B^+(R_1/2))\times ]-3/2,0[$, $x_{03}=0$.
 From the $\varepsilon$-regularity theory developed in \cite{CKN} and \cite{S3}, \cite{S4}, see details also in \cite{ESS2003},
in particular, we can show the validity of (\ref{maxu})
in $(\mathbb R^3_+\setminus B^+(3R_1/2))\times ]-4/3,0[$.

The second statement of the proposition can be deduced  from the local regularity theory for the heat equation in the following sense. This follows from bootstrap arguments involving the vorticity equation and is described in detail in the Lemma \ref{vortbootstrap} found in the Appendix. $\Box$

\setcounter{equation}{0}
\section{Rescaling. Scenario I}
Let us go back to our original problem (\ref{directsystem})--(\ref{directic}).

We assume that $T>0$ is a blowup time. Theorem \ref{limhalf} can be proven ad absurdum. Suppose that there exists a sequence $t_n\uparrow T$ such that 
$$M:=\sup\limits_{n}\|v(\cdot,t_n)\|_{3,\mathbb R^3_+}<\infty,$$
then $T$ is NOT a blowup time.

It is known that there exists a global weak Leray-Hopf solution (energy solution) to initial boundary value problem (\ref{directsystem})--(\ref{directic}). This solution coincides with $v$ on the interval $]0,T[$ and that is why we are going to denote it still by $v$. 
 Arguments similar to used in the previous section show that for every $\epsilon>0$ there exists $R_1(\epsilon)>0$ such that $$\sup\limits_{x\in \mathbb R^3_+\setminus B^+(R_1), \epsilon\leq t\leq T}|v(x,t)|<\infty.$$ So, by the definition of blowup time $T$,
 there should a singular point $x_0=(x_0,x_{03})\in \mathbb R^3_+$ at $t=T$, i.e., a point such that $v\notin  L_\infty(B(x_0,r)\cap\mathbb R^3_+)\times ]T-r^2,T[)$ for any positive $r$. Without loss of generality, we may assume that $x_0'=0$. Then one should consider two case
$$x_{03}=0$$
and
$$x_{03}>0.$$

Now, let us focus on the first case. We know from \cite{S3} and \cite{S4}   that it must be
$$\frac 1{a^2}\int\limits_{Q^+(a)}(|v|^3+|q-[q]_{B^+(a)}|^\frac 32)dx dt>\varepsilon$$
for all $0<a<a_0$, for some positive $a_0$, and for some universal constant
$\varepsilon$.

In this scenario, our rescaling  will be as follows:
$$u^{(n)}(y,s)=\lambda_nv(x,t), \qquad
p^{(n)}(y,s)=\lambda^2_nq(x,t),$$
where
$$x=\lambda_ny,\qquad t=T+\lambda^2_ns,\qquad \lambda_n=\sqrt{\frac {T-t_n}2}$$

So, sufficiently smooth solutions $u^{(n)}$ and $p^{(n)}$  satisfy  (\ref{nsystem})-(\ref{nic}) with   
$$u^{(n)}_0(y)=\lambda_nv(\lambda_ny,t_n).$$
Hence, $$\sup\limits_n\|u^{(n)}_0\|_{3,\mathbb R^3_+}=
M.$$
Without loss of generality, we may assume
$$u^{(n)}_0\rightharpoonup u_0$$
in $L_3(\mathbb R^3_+)$.
So, all assumptions and statements  of Proposition \ref{convergence} and Proposition \ref{decaypro}
hold for $u^{(n)}$ and $p^{(n)}$ and for their limits $u$ and $p$. 

Now, we shall show that
\begin{equation}\label{zeroattheend}
u(x,0)=0\end{equation}
 Indeed, it is not difficult to see that
$$\frac 1{a^\frac {15}8}\int\limits_{B^+(a)}|u^{(n)}(x,0)|^\frac 98dx\to \frac 1{a^\frac {15}8}\int\limits_{B^+(a)}|u(x,0)|^\frac 98dx.$$
And, on the other hand,
$$\frac 1{a^\frac {15}8}\int\limits_{B^+(a)}|u^{(n)}(x,0)|^\frac 98dx\leq \Big(\int\limits_{B^+(a)}|u^{(n)}(x,0)|^3dx\Big)^\frac 38=$$$$=\Big(\int\limits_{B^+(\lambda_na)}|v(x,0)|^3dx\Big)^\frac 38\to0$$
as $n\to\infty$. The latter is true as we can easily show that the integral
$$\int\limits_{\mathbb R^3_+}|v(x,0)|^3dx$$
is finite.

Now, we need to show that the limit function is not identically zero.
Fix $0<a_*<1/4$ 
then we have
$$\frac 1{a^2}\int\limits_{Q^+(a)}(|u^{(n)}|^3+|p^{(n)}-[p^{(n)}]_{B^+(a)}|^\frac 32)dxdt=$$
$$=\frac 1{(\lambda_na^2)}\int\limits_{Q^+(\lambda_na)}(|v|^3+|q-[q]_{B^+(\lambda_na)}|^\frac 32)dxdt>\varepsilon$$
for $0<a<a_*$.

We know that
\begin{equation}\label{unsup}
M_1(a_*):=\sup\limits_n\Big\{\frac 1{(2a)^2_*}\int\limits_{Q^+(2a_*)}(|u^{(n)}|^3+|p^{(n)}-[p^{(n)}]_{B^+(2a_*)}|^\frac 32)dxdt+
$$$$+\sup\limits_{-(2a_*)^2<t<0}\|u^{(n)}(\cdot,t)\|^2_{2,B^+(2a_*)}+\|\nabla u^{(n)}\|^2_{2,Q^+(2a_*)}\Big\}<\infty.
\end{equation}

 Let us fix a $C^2$-domain $\Omega_*$ such that $B^+(a_*)\subset\Omega_*\subset B^+(2a_*)$ and let  $Q_*=\Omega_*\times ]-a_*^2,0[$. We may use the same type of decompositions as in the previous sections
$$u^{(n)}=w^1+w^2,\qquad p^{(n)}=r^1+r^2$$
so that
$$\partial_tw^1-\Delta w^1+\nabla r^1=0,\qquad \mbox{div}\,w^1=0$$
in $Q_*$,
$$w^1(x',0,t)=0$$
for all $(x,t)\in\partial\Omega_*\times [-a_*^2,0]$ with $x_3=0$, 
$$w^1(x,-a_*^2)=u^{(n)}(x,-a_*^2)$$
for all $x\in \Omega_*$ and
$$\partial_tw^2-\Delta w^2+\nabla r^2=-\mbox{div}\,u^{(n)}\otimes u^{(n)},\qquad \mbox{div}\,w^2=0$$
in $Q_*$,
$$w^2=0$$
on the parabolic boundary of $Q_*$.


By the Solonnikov coercive estimate, see \cite{Sol1973} and \cite{Sol2003UMN}, we have
\begin{equation}\label{soloncoer}
\|w^{2}\|_{W^{2,1}_{\frac{12}{11},\frac{3}{2}}(Q_{*})}+\|r^{2}\|_{W^{1,0}_{\frac{12}{11},\frac{3}{2}}(Q_{*})}\leqslant C(a_{*})
\|u^{(n)}\cdot \nabla u^{(n)}\|_{L_{\frac{12}{11},\frac{3}{2}}(Q_{*})}.
\end{equation}
Using (\ref{soloncoer}) one infers
$$\int\limits_{Q^+(a_*)}|r^2-[r^2]_{B^+(a_*)}|^\frac 32 dx dt\leq 
 c(a_*)
\int\limits^0_{-a_*^2}\Big(\int\limits_{B^+(a_*)}|\nabla r^2|^\frac {12}{11}dx\Big)^\frac {11}8dt\leq $$$$\leq c(a_*)
\int\limits^0_{-a_*^2}\Big(\int\limits_{\Omega_*}|\nabla r^2|^\frac {12}{11}dx\Big)^\frac {11}8dt\leq $$
$$\leq  c(a_*)
\int\limits^0_{-a_*^2}\Big(\int\limits_{\Omega_*}|u^{(n)}\cdot \nabla u^{(n)}|^\frac {12}{11}dx\Big)^\frac {11}8dt\leq $$
$$\leq c(a_*)\|\nabla u^{(n)}\|^\frac 32_{2,Q_*}\| u^{(n)}\|^\frac{3}{4}_{2,\infty,Q_*}\| u^{(n)}\|^\frac 34_{3,Q_*}$$

\begin{equation}\label{presest}
\leq c(a_*,M_1)\| u^{(n)}\|^\frac 34_{3,Q_*}\leq c(a_*,M_1)\Big(
\int\limits_{Q^+(2a_*)}|u^{(n)}|^3dx dt\Big)^\frac 14.
\end{equation}
By the local regularity theory up the boundary for the Stokes system  developed in \cite{S3} and \cite{S4} and by (\ref{unsup})-(\ref{presest}), 
 for any $s>\frac {12}{11}$,
$$\int\limits^0_{-(a_*/2)^2}\Big(\int\limits_{B^+(a_*/2)}|\nabla r^1|^sdx\Big)^\frac 3{2s}dt\leq $$$$\leq c(a_*,s)(\|w^{1}\|^{\frac 32}_{L_{\frac{12}{11},\frac 32}(Q^+(a_*))}+\|\nabla w^{1}\|^{\frac 32}_{L_{\frac{12}{11},\frac 32}(Q^+(a_*))}+$$$$+\|r^{1}-[r^1]_{B^*(a_*)}\|^{\frac 32}_{L_{\frac{12}{11},\frac 32}(Q^+(a_*))})\leq  $$
$$\leq c(a_*,s)(\|w^{2}\|^{\frac 32}_{L_{\frac{12}{11},\frac 32}(Q^+(a_*))}+\|\nabla w^{2}\|^{\frac 32}_{L_{\frac{12}{11},\frac 32}(Q^+(a_*))}+$$$$+\|u^{n}\|^{\frac 32}_{L_{\frac{12}{11},\frac 32}(Q^+(a_*))}+\|\nabla u^{n}\|^{\frac 32}_{L_{\frac{12}{11},\frac 32}(Q^+(a_*))}+$$$$+\int\limits_{Q^+(a_*)}|r^1-[r^1]_{B^*(a_*)}|^\frac 32 dx dt)\leq  $$
$$\leq c(a_*,M_{1},s)(\int\limits_{Q^+(a_*)}(|p^{(n)}-[p^{(n)}]_{B^+(a_*)}|^\frac 32+|r^2-[r^2]_{B^*(a_*)}|^\frac 32)dx dt+1)\leq $$$$ \leq c(a_*, M_1).$$

Next, we have, for any $0<a<a_*/2$ and for any $s=9$,

$$\varepsilon<\frac 1{a^2}\int\limits_{Q^+(a)}(|p^{(n)}-[p^{(n)}]_{B^+(a)}|^\frac 32+|u^{(n)}|^3)dxdt\leq $$$$\leq
\frac 1{a^2}\int\limits_{Q^+(a)}(|r^1-[r^1]_{B^+(a)}|^\frac 32)dxdt+ $$
$$+
\frac 1{a^2}\int\limits_{Q^+(a)}(|r^2-[r^2]_{B^+(a)}|^\frac 32+|u^{(n)}|^3)dxdt\leq $$$$\leq c(a_*,M_1)a^2+
\frac 1{a^2}\int\limits_{Q^+(a)}(|r^2-[r^2]_{B^+(a)}|^\frac 32+|u^{(n)}|^3)dxdt\leq $$
$$\leq \frac {c(a_*)}{a^2}\int\limits_{Q^+(a_*)}(|u^{(n)}|^3+|r^2-[r^2]_{B^+(a_*)}|^\frac 32)dxdt+c(a_*,M_1)a^2\leq $$
$$\leq \frac {c(a_*)}{a^2}\int\limits_{Q^+(a_*)}|u^{(n)}|^3dx dt+\frac {c(a_*,M_1)}{a^2}\Big(
\int\limits_{Q^+(2a_*)}|u^{(n)}|^3dx dt\Big)^\frac 14+c(a_*,M_1)a^2\leq $$
$$\leq\frac {c(a_*,M_1)}{a^2}\Big(
\int\limits_{Q^+(2a_*)}|u^{(n)}|^3dx dt\Big)^\frac 14+c(a_*,M_1)a^2.$$
For sufficiently small $a>0$,
$$0<\varepsilon/2<\varepsilon-c(a_*,M_1)a^2\leq \frac {c(a_*,M_1)}{a^2}\Big(
\int\limits_{Q^+(2a_*)}|u^{(n)}|^3dx dt\Big)^\frac 14$$
and thus
$$\int\limits_{Q^+(2a_*)}|u^{(n)}|^3dx dt>\Big(\frac {a^2
\varepsilon}{c(a_*,M_1)}\Big)^4.$$
Passing to the limit as $n\to\infty$, we find
\begin{equation}\label{non-triviality}\int\limits_{Q^+(2a_*)}|u|^3dx dt>\Big(\frac {a^2
\varepsilon}{c(a_*,M_1)}\Big)^4.\end{equation}

Next, we follow arguments of the paper \cite{ESS2003} that related with backward uniqueness for the heat operator with lower order terms.
Indeed, by Proposition \ref{decaypro}, we have
$$|\partial_t\omega-\Delta \omega|\leq c(|\omega|+|\nabla\omega|)$$
in $\{x\in\mathbb R^3:\,\,x_3>2R_1\}\times ]-5/4,0[$, where $\omega=\nabla \wedge u$. Then, because of (\ref{zeroattheend}), we can state that
\begin{equation}\label{vorticityzero}\omega(x,t)=0\end{equation}
for all $(x,t)\in \{x\in\mathbb R^3:\,\,x_3>2R_1\}\times ]-6/5,0[$. 

Applying unique continuation through spatial boundaries, we may
conclude that (\ref{vorticityzero}) is valid in $(\mathbb R^3_+\setminus B^+(R_1)\times ]-7/6,0[$. 

For the final component of the proof, we initially refer back to the section on a priori estimates. Indeed, we use the same decomposition applied to the scaled solutions:
$$u^{(n)}:= u^{1,(n)}+u^{2,(n)}.$$ The Solonnikov estimates implies the following inequalities for $k=0,1\ldots$
\begin{equation}\label{solonspacederbdd}
\|\nabla^{k}u^{1,(n)}(\cdot,t)\|_{\infty,\mathbb{R}^{3}_{+}}\leq\frac{c}{{(t+2)}^{\frac{k+1}{2}}}\|u_{0}^{(n)}(\cdot)\|_{3,\mathbb{R}^{3}_{+}}\leq\frac{cM}{{(t+2)}^{\frac{k+1}{2}}},
\end{equation}
\begin{equation}\label{solonspacederL3}
\|\nabla^{k}u^{1,(n)}(\cdot,t)\|_{3,\mathbb{R}^{3}_{+}}\leq\frac{c}{{(t+2)}^{\frac{k}{2}}}\|u_{0}^{(n)}(\cdot)\|_{3,\mathbb{R}^{3}_{+}}\leq\frac{cM}{{(t+2)}^{\frac{k}{2}}},
\end{equation}

for any $t\in ]-2,0[.$ One can observe that $u^{2,(n)}$ satisfies the following:
\begin{equation}\label{u_{2}eqn}
\partial_tu^{2,(n)}+\mbox{div}\,u^{2,(n)}\otimes u^{2,(n)}-\Delta u^{2,(n)}=-\nabla p^{2,(n)}+f^{2,(n)},\qquad
\end{equation}
$$ \mbox{div}\,u^{2,(n)}=0$$
in $Q^+_{-2,0}$, the boundary conditions 
$$u^{2,(n)}(x',0,t)=0$$
for $(x',t)\in\mathbb R^2\times [-2,0]$, and the initial conditions
$$u^{2,(n)}(\cdot,-2)=0$$
in $R^3_+$.
Here,
\begin{equation}
f^{2,(n)}:=\mbox{div}\,(u^{1,(n)}\otimes u^{1,(n)}+u^{2,(n)}\otimes u^{1,(n)}+u^{1,(n)}\otimes u^{2,(n)}).
\end{equation}
Using (\ref{solonspacederbdd})-(\ref{solonspacederL3}) together with interpolation and the apriori estimates previously obtained for $u^{2,(n)}$, it is not so difficult to see 
$$\|f^{2,(n)}\|_{2,Q_{-\frac{3}{2},0}^{+}}\leqslant C(M).$$
For $N>0$ define $$Q_{-\frac{3}{2},0}^{+}(N):= B^{+}(N)\times ]-\frac{3}{2},0[.$$
Using standard arguments, we claim (up to subsequence):
\begin{equation}\label{L2,infntyweakstar}
u^{2,(n)}\stackrel{*}{\rightharpoonup} u^{2}
\end{equation}
in $L_{2,\infty}(Q_{-\frac{3}{2},0}^{+})$,
\begin{equation}\label{weakconverggrad}
\nabla u^{2,(n)}\rightharpoonup\nabla u^{2}
\end{equation}
in $L_{2}(Q_{-\frac{3}{2},0}^{+})$,
\begin{equation}\label{presweakconverg}
p^{2,(n)}-[p^{2,(n)}]_{B^{+}(N)}\rightharpoonup p_{N}
\end{equation}
in $L_{\frac{3}{2}}(Q_{-\frac{3}{2},0}^{+}(N))$ (for $N=1,2\ldots$),
\begin{equation}\label{preslimitest}
\|p_{N}\|_{L_{Q_{-\frac{3}{2},0}^{+}(N)}}\leqslant C(M)(N^{\frac{1}{3}}+N+N^{\frac{1}{2}}),
\end{equation}
\begin{equation}\label{L3strongconverg}
 u^{2,(n)}\rightarrow u^{2}
\end{equation}
in $L_{3}(Q_{-\frac{3}{2},0}^{+}(N))$ (for $N=1,2\ldots$),
\begin{equation}\label{weakL2pointwise}
u^{2,(n)}(\cdot,t)\rightharpoonup u^{2}(\cdot,t)
\end{equation}
in $L_{2}(\mathbb{R}^{+}_{3})$ for each $t\in]-\frac{3}{2},0[$,
\begin{equation}\label{weakconvergf}
 f^{2,(n)}\rightharpoonup f^{2}
\end{equation}
in $L_{2}(Q_{-\frac{3}{2},0}^{+})$.
Upon passage to the limit we also obtain that (after appropriate adjustment of $u^{2}(\cdot,t)$ on a subset of $[-\frac{3}{2},0]$ of Lebesgue measure zero) that for any $w\in L_{2}(\mathbb{R}^{3}_{+})$ the function
\begin{equation}\label{weakL2contu2}
t:\rightarrow \int\limits_{\mathbb{R}^{3}_{+}}u^{2}(x,t)w(x)dx
\end{equation}
is in $C([-\frac{3}{2},0])$.
Now, we see from (\ref{solonspacederbdd}) that for $k=0,1,2,\ldots$ (up to subsequence):
\begin{equation}\label{u1weakstarconverg}
\nabla^{k}u^{1,(n)}\stackrel{*}{\rightharpoonup}\nabla^{k}u^{1},
\end{equation}
 in $L_{\infty}(Q_{-\frac{3}{2},0}^{+})$ along with estimate
\begin{equation}\label{u1spaceder}
\|\nabla^{k}u^{1}\|_{L_{\infty}(Q_{-\frac{3}{2},0}^{+})}\leqslant c_{k}M
\end{equation}
for universal constants $c_{k}$. It can be seen that weak star convergence occurs in  $L_{3,\infty}(Q_{-\frac{3}{2},0}^{+})$ with analogous estimates. Upon passage to the limit we also obtain that (after appropriate adjustment of $u^{2}(\cdot,t)$ on a subset of $[-\frac{3}{2},0]$ of Lebesgue measure zero) that for any $w\in L_{\frac{3}{2}}(\mathbb{R}^{3}_{+})$ the function
\begin{equation}\label{weakL3contu1}
t:\rightarrow \int\limits_{\mathbb{R}^{3}_{+}}u^{2}(x,t)w(x)dx
\end{equation}
is in $C([-\frac{3}{2},0])$.
We obtain that the limit functions $u$ and $f^{2}$  can be decomposed, for $(x,t)\in Q_{-\frac{3}{2},0}^{+}$ as  follows
\begin{equation}\label{udecom}
u(x,t):= u^{1}(x,t)+u^{2}(x,t),
\end{equation}
\begin{equation}\label{f2decomp}
f^{2}:=\mbox{div}\,(u^{2}\otimes u^{1}+u^{1}\otimes u^{2})+F^{2}_{1}.
\end{equation}
Furthermore, the following estimate for $F^{2}_{1}$ is valid for $2\leqslant p\leqslant \infty$ and $k=0,1\ldots$
\begin{equation}\label{F2spacederest}
\|\nabla^{k}F^{2}_{1}\|_{L_{p}(Q_{-\frac{3}{2},0}^{+})}\leqslant c(k,p)M^{2}.
\end{equation}
Furthermore, from (\ref{weakL2contu2}) and (\ref{weakL3contu1}), we observe that (after appropriate adjustment of $u(\cdot,t)$ on a subset of $[-\frac{3}{2},0]$ of Lebesgue measure zero) that for any $\phi\in C^{\infty}_{0}(\mathbb{R}^{3})$ the function
\begin{equation}\label{weakcontu}
t:\rightarrow \int\limits_{\mathbb{R}^{3}_{+}}u^{2}(x,t)\phi(x)dx
\end{equation}
is in $C([-\frac{3}{2},0])$.\\
 Using smoothness properties of $u^{2,(n)}$, along with  convergence facts and properties of limit functions (described in (\ref{L2,infntyweakstar})-(\ref{weakL2contu2})), we claim there exists a set $\Sigma\subset ]-\frac{3}{2},0[$ of full measure, i.e., $|\Sigma|=\frac{3}{2}$, such that for $t_{0}\in\Sigma$:
 $$\|u^{2}(\cdot,t_{0})\|_{L_{2}(\mathbb{R}^{3}_{+})},\|\nabla u^{2}(\cdot,t_{0})\|_{L_{2}(\mathbb{R}^{3}_{+})}<\infty.$$
 Moreover $u^{2}$ is a weak Leray-Hopf solution to the following initial value problem on $\mathbb{R}^{3}_{+}\times]t_{0},0[$: 
 \begin{equation}\label{u_{2}eqn}
\partial_tU+\mbox{div}\,U\otimes U-\Delta U=-\nabla P+f^{2},\qquad
\end{equation}
$$ \mbox{div}\,U=0,$$
$$U(x',0,t)=0$$
and 
$$U(\cdot,t_{0})= u^{2}(\cdot,t_{0}).$$

 The initial value and source are $u(\cdot,t_{0})$ and $f^{2}$ respectively.
 Then, by the short time unique solvability results for the Navier-Stokes system in unbounded domains with smooth boundary (see \cite{Heywood1980} and \cite{L1970}, for example), we can find a number $\delta_{0}>0$ such that
 $$\partial_{t}u^{2},\, \nabla^{2}u^{2},\,\nabla p_{N}\in L_{2}(\mathbb{R}^{3}_{+}\times ]t_{0},t_{0}+\delta_{0}[).$$
 Consequently one may use the parabolic embedding theorems, together with the regularity theory for linear systems and bootstrap arguments, to obtain (for arbitrary $\epsilon>0$):
 $$\sup_{t_{0}+\epsilon<t<t_{0}+\delta_{0}}\sup_{\mathbb{R}^{3}_{+}}|u^{2}(x,t)|\leqslant c.$$
 Using properties of $u^{1}$ and  similar arguments to Lemma \ref{vortbootstrap} from the Appendix, we obtain for any $\delta>0$, $k=0,1,\ldots$
 $$\sup_{t_{0}+2\epsilon<t<t_{0}+\delta_{0}}\sup_{\mathbb{R}^{3}_{+\delta}}|\nabla^{k}u(x,t)|\leqslant c_{1}(\delta,k,\epsilon,c,\|\nabla u\|_{L_{2,unif}(Q_{-\frac{3}{2},0}^{+})}).$$
Additionally, we obtain sufficient regularity on the time derivative of the vorticity to apply the unique continuation theorem through spatial boundaries. Repeating arguments in \cite{ESS2003}, obtain $\omega(\cdot,t)=0$ in $\mathbb{R}^{3}_{+\delta}\times]t_{0}+2\epsilon,t_{0}+\delta_{0}[$ and arbitrary $\delta>0$.
Hence, for a.a $t\in]t_{0}+2\epsilon,t_{0}+\delta_{0}[$, $u$ is a harmonic function, which satisfies the boundary condition $u(x,t)=0$ if $x_{3}=0$. But for a.a $t\in]t_{0}+2\epsilon,t_{0}+\delta_{0}[$, $L_{3}$-norm of $u$ over $\mathbb{R}^{3}_{+}$ is finite. This leads to the conclusion that, for the same $t$, $u(\cdot,t)=0$ in $\mathbb{R}^{3}_{+}$.
Exploiting the continuity of $u$ described in (\ref{weakcontu}) and arbitrariness of $\epsilon$, it is simple to see that $u(\cdot,t)=0$ for all $t\in [t_{0},t_{0}+\delta_{0}].$ Since $t_{0}\in |\Sigma$ with $|\Sigma|=\frac{3}{2}$, it is immediate that one obtaines the same conclusion for every $t\in]-\frac{3}{2},0[$.
    
 The latter contradicts with (\ref{non-triviality}) and thus $T$ is not a blowup time.$\Box$

\setcounter{equation}{0}
\section{Rescaling. Scenario II}
Here, the scaling is $x=x_0+\lambda_n y$. So, we replace $\mathbb R^3_+$ with $\mathbb R^3_{h}=\{y=(y',y_3)\in\mathbb R:\,\,y_3>h\}$ with $h=h_n=-x_{03}/\lambda_n$.

 In the case, sufficiently smooth functions $u^{(n)}$ and $p^{(n)}$  are a solution to   the following initial boundary value problem:
$$\partial_tu^{(n)}+\mbox{div}\,^{(n)}\otimes u^{(n)}-\Delta u^{(n)}=-\nabla p^{(n)},\qquad \mbox{div}\,u=0$$
in $\mathbb R^3_{h_n}\times ]-2,0[$, 
$$u^{(n)}(x',-h_n,t)=0$$
for $(x',t)\in\mathbb R^2\times [-2,0]$,
$$u^{(n)}(\cdot,-2)=u^{(n)}_0(\cdot)\in L_3(\mathbb R^3_{h_n})$$
and
$$\sup\limits_n\|u^{(n)}_0\|_{3,\mathbb R^3_{h_n}}\leq
M.$$
Without loss of generality, we may assume
$$u^{(n)}_0\rightharpoonup u_0\in L_3(\mathbb R^3)$$
in $L_3(\mathbb R^3_h)$ for any $h>-\infty$.

 We can use estimates of Section 2 in domains $\mathbb R^3_{h_n}$ with constants independent of $n$.
\begin{pro}\label{convergence2}
There exist subsequences still denoted in the same way with the following properties:
\begin{equation}\label{weakconvergence2}
u^{(n)}\rightharpoonup u
\end{equation}
in $L_\frac {10}3(R^3_h\times ]-2,0[)$ for any $h>-\infty$ with $u\in L_\frac {10}3(Q_{-2,0})$ and $Q_{-2,0}=\mathbb R^3\times ]-2,0[$,
\begin{equation}\label{weakconvergencenabla2}\nabla u^{(n)}\rightharpoonup \nabla u\end{equation}
in $L_2(B(R)\times ]-2+\delta,0[)$ for any $R>0$ and any $0<\delta<2$,
\begin{equation}\label{strongconvergence2}
u^{(n)}\rightarrow u
\end{equation}
in $L_3(B(R)\times ]-3/2,0[)$ for any $R>0$;
\begin{equation}\label{weakconverpressure2}
p^{(n)}\rightharpoonup p
\end{equation}
in $L_{\frac 32}(B(R)\times ]-3/2,0[)$ for any $R>0$.

Functions $u$ and $p$ satisfy the Navier-Stokes system in $\mathbb R^3\times ]-3/2,0[$. 

For the pressure $p$, the following global estimates are valid:
$$p=p^1+p^2$$
and
$$p^{2}=\sum\limits^5_{i=1}p^{2,i}.
$$
with the estimates
$$\|\nabla p^1\|_{3,\mathbb R^3\times ]-3/2,0[}+\|\nabla p^{2,1}\|_{2,\mathbb R^3\times ]-3/2,0[}+$$
\begin{equation}\label{gradpressure2}
+\|\nabla p^{2,2}\|_{9/8,3/2,\mathbb R^3\times ]-3/2,0[}
+\|\nabla p^{2,3}\|_{3/2,\mathbb R^3\times ]-3/2,0[}+\end{equation}
$$
+\|\nabla p^{2,4}\|_{6/5,3/2,\mathbb R^3\times ]-3/2,0[}
+\|\nabla p^{2,5}\|_{3/2,\mathbb R^3\times ]-3/2,0[}<\infty.$$

Moreover, for any $R>0$, the limits pair $u$ and $p$ satisfies the local energy inequality
$$\int\limits_{B(x_0,R)}|\varphi^2(x,t)|u(x,t)|^2dx +2\int\limits^{t}_{t_0-R^2}\int\limits_{B(x_0,R)}\varphi^2|\nabla u|^2dxdt   \leq $$
\begin{equation}\label{locenergy2}
\leq\int\limits^{t}_{t_0-R^2}\int\limits_{B(x_0,R))}\Big(|u|^2(\partial_t\varphi^2+\Delta \varphi^2u^\cdot \varphi^2(|u|^2+2p)\Big) dxds
\end{equation}
for all $-3/2<t_0-R^2<t\leq t_0\leq0$, for all $x_0\in R^3$, and for all $\varphi\in C^\infty_0(B(x_0,R)\times ]t_0-R^2,t_0+R^2[)$.\end{pro}
 The proof of Proposition \ref{convergence} goes along the lines of the proof of Proposition \ref{convergence} with minor modifications.

A major simplification in Scenario II, compared with Scenario I, is related to showing non-triviality of the limit solution. In the interior case we may follow the local pressure decomposition used in \cite{Ser2012} one of which is harmonic and the other satisfies a coercive estimate. Interior properties of harmonic functions are essential in the use of this decomposition. In the boundary case of Scenario II the same decomposition doesn't apply. Instead, one uses a local decomposition of the velocity and pressure together with estimates for the Stokes system near the boundary as described in Scenario I. The remainder of the proof is similar to that described for Scenario I, with few minor modifications. $\Box$ 

\section{Appendix }
The following Lemma seems to be known. It is useful for verifying the hypothesis for the theorems of backward uniqueness and unique continuation through spatial boundaries of parabolic operators. We give a proof for the readers convenience.
\begin{lemma}\label{vortbootstrap}
Let $u$ belong to $L_{2}(B^{+}(R)\times ]-\frac{3}{2},0[)$   for any $R>0$. Also let 
\begin{equation}\label{gradunif}
\|\nabla u\|_{L_{2,unif}(Q^{+}_{-\frac{3}{2},0})}:=\sup_{x_{0}\in\mathbb{R}^{3}_{+}}\|\nabla u\|_{L_{2}(B(x_{0},1)\times ]-\frac{3}{2},0[)}<\infty.
\end{equation}
Suppose that  Functions $u$ and $p$ satisfy (\ref{system}) in $\mathbb R^3_+\times ]-3/2,0[$.
Furthermore, suppose that
\begin{equation}\label{ubddaway}
|u(x,t)|\leq c
\end{equation}
for all $(x,t)\in (\mathbb R^3_+\setminus B^+(R_1))\times ]-5/4,0[$ and for some universal constant $c$.
Then we infer that, given any $\delta>0$ and $k= 1,2,\ldots$, there exists a constant $c_{1}(\delta,k,c,\|\nabla u\|_{L_{2,unif}(Q^{+}_{-\frac{3}{2},0})})>0$ such that 
\begin{equation}\label{uspacialsmooth}
|\nabla^{k}u(x,t)|\leqslant c_{1}(\delta,k,c,\|\nabla u\|_{L_{2,unif}(Q^{+}_{-\frac{3}{2},0})})
\end{equation}
for all $(x,t)\in(\mathbb{R}^{3}_{+\delta}\setminus B^{+}(2R_{1})\times]-\frac{5}{4},0[$.
\end{lemma}
\textsl{Proof of Lemma \ref{vortbootstrap}}\
Let $\omega$ denote the vorticity, namely $\omega:=\nabla \wedge u$.
It satisfies in $\mathbb{R}^{3}\times]-\frac{3}{2},0[$:
\begin{equation}\label{vorteqn1}
\partial_t\omega-\Delta \omega=\mbox{div}(\omega\otimes u-u\otimes \omega),
 \end{equation}
 and
 \begin{equation}\label{vorteqn2}
 -\Delta u:=\nabla\wedge \omega.
 \end{equation}
 Let $x_{0}\in\mathbb{R}^{3}_{+\frac{\delta}{2}}\setminus B^{+}(\frac{3R_{1}}{2})$ and let $a_{1}$ be sufficiently small such that $B(x_{0},a)\in\mathbb{R}^{3}_{+{\delta}}\setminus B^{+}(2R_{1}) $.
 Using (\ref{gradunif}) and local regularity theory for heat equation (e.g Appendix of \cite{NRS1996}) and a parabolic embedding theorem found in \cite{LSU}, obtain (for $a_{1}<a$ and $\tau_{1}>\frac{3}{2}$):
\begin{equation}\label{estgradvort}
\|\nabla\omega\|_{2,B(x_{0},a_{1})\times ]-\tau_{1},0[)}\leqslant c_{2}(c_{1},a,a_{1},\tau_{1}),
\end{equation}
\begin{equation}\label{vortparaembed1}
\|\omega\|_{\frac{10}{3},B(x_{0},a_{1})\times ]-\tau_{1},0[)}\leqslant c_{2}(c_{1},a,a_{1},\tau_{1}).
\end{equation}
Using these estimates along with \ref{vorteqn2} and local regularity for Laplace equation obtain:
\begin{equation}\label{estu1}
\|\nabla^{2}u\|_{2,B(x_{0},a_{1})\times ]-\tau_{1},0[)}+\|\nabla u\|_{\frac{10}{3},B(x_{0},a_{1})\times ]-\tau_{1},0[)}\leqslant c_{3}(c_{1},a,a_{1},\tau_{1}).
\end{equation}
And thus
$$\mbox{div}(\omega\otimes u-u\otimes \omega)\in L_\frac 53(B(x_0,a_1)\times ]-\tau_1,0[).$$
Then the local regularity theory gives
$$\|\partial_t\omega\|_{5/3,B(x_0,a_2)\times ]-\tau_2,0[}+\|\nabla^2\omega\|_{5/3,B(x_0,a_2)\times ]-\tau_2,0[}\leq c_{4}(c_{1},a,a_1,a_2,\tau_1,\tau_2)$$
for any $0<a_2<a_1$ and $0<\tau_2<\tau_1$. Now, according to the parabolic embedding theorem, see
\cite{LSU}, we have
$$\omega\in L_5(B(x_0,a_2)\times ]-\tau_2,0[), \qquad \nabla\omega\in L_\frac 52(B(x_0,a_2)\times ]-\tau_2,0[).$$
with the corresponding estimates. The same estimates are valid for $\nabla u$ instead of $\omega$.
So,
$$\|\partial_t\omega\|_{5/2,B(x_0,a_3)\times ]-\tau_3,0[}+\|\nabla^2\omega\|_{5/2,B(x_0,a_3)\times ]-\tau_3,0[}\leq c(c_{1},a,a_1,a_2,a_3,\tau_1,\tau_2,\tau_3)$$
for any $0<a_3<a_2$ and $0<\tau_3<\tau_2$. From the parabolic embedding theorem , we find that $\omega\in L_s(B(x_0,a_3)\times ]-\tau_3,0[)$ with any $s>1$ and thus $$\omega\otimes u-u\otimes\omega\in L_s(B(x_0,a_3)\times ]-\tau_3,0[)$$ with any $s>1$. Taking $s>5$, $a_{4}<a_{3}$ and $\tau_{4}<\tau_{3}$, we can use the local regularity theory for the heat equation and embeddings once more to obtain that $\nabla\omega\in L_s(B(x_0,a_4)\times ]-\tau_4,0[)$ with any $s>1$ and $\omega$ is H\"older continuous in the same domain with the required estimates. Thus, we infer
$$\mbox{div}(\omega\otimes u-u\otimes \omega)\in L_s(B(x_0,a_4)\times ]-\tau_4,0[),$$
for any $s>1$. Applying local regularity theory for the heat equation one more time, we get $$\omega\in W^{2,1}_{s}(B(x_{0},a_{5})\times]-\tau_{5},0[)$$ with required constant dependence. Here, $s>1$ is arbitrary and $a_{5},\,\tau_{5}<a_{4},\,\tau_{4}$. Using the parabolic embedding theorem one more time gives that for $s<5$ we in fact have that $\nabla\omega$ is H\"older continuous in the same domain. Now for $a_{6}<a_{5}$, (\ref{vorteqn2}) together with local regularity for the Laplace equation gives $\ref{uspacialsmooth}$ for $k=1$ and the same estimate for $\|\nabla^{2}u\|_{L_{s,\infty}(B(x_{0},a_{5})\times]\-tau_{5},0[)}$ ($s>1$ is arbitrary). These conclusions easily allow us to iterate the same arguments to spacial derivatives of any order. $\Box$

\end{document}